\documentclass[11pt,a4paper]{article}

\usepackage[T1]{fontenc}
\usepackage[utf8]{inputenc}
\usepackage{lmodern}
\usepackage[english]{babel}
\usepackage{amsmath,amssymb}
\usepackage{textalpha} % Greek words in running text
\usepackage{csquotes}
\usepackage[margin=2.7cm]{geometry}

\usepackage[backend=biber,
            style=authoryear,
            natbib=true,
            maxcitenames=2,
            uniquename=false,
            uniquelist=false,
            dashed=false]{biblatex}
\usepackage[hidelinks]{hyperref}
\title{A Rightly Forgotten Mathematician?\\[0.3em]
\large József Suták (1865--1954) and Budapest Mathematics\\ between Bolyai and the Numerus Clausus}
\author{András Bátkai\\[0.3em]
\small Pädagogische Hochschule Vorarlberg, Liechtensteiner Straße 33--37, 6800 Feldkirch, Austria\\
\small \texttt{andras.batkai@ph-vorarlberg.ac.at}}
\date{}

\begin{document}

\maketitle

\begin{abstract}
József Suták (1865--1954), a Piarist priest, \emph{gymnasium} teacher and from 1912 full professor of geometry at the University of Budapest, is scarcely more than a footnote in the history of Hungarian mathematics. His name survives mainly because he co-signed, with Lipót Fejér, the acceptance of John von Neumann's 1926 dissertation on the axiomatic construction of set theory. Drawing on his printed works, on the records of the Budapest Faculty of Philosophy and of the Catholic Saint Stephen Academy, and on a hitherto unknown part of his personal papers, this article reconstructs his career. It follows his Hungarian translation of Bolyai's \emph{Appendix} (1897), his empiricist lecture on the axioms of geometry (1898), his papers on the logical problems of set theory (1911) and on the axiomatic construction of projective geometry (1914), and his later differential-geometric work. It analyses his appointment in 1911/12, when the faculty preferred him to Frigyes Riesz and Alfréd Haar, and his role as chairman of the admissions committee in implementing the numerus clausus after 1920. The answer to the title's question is twofold: as a mathematician, Suták is largely rightly forgotten; for the history of Budapest mathematics and its university, he is not.

\medskip
\noindent\textbf{Keywords:} József Suták; János Bolyai; foundations of mathematics; set theory; projective geometry; numerus clausus; University of Budapest; John von Neumann; Piarist order

\end{abstract}

\section{Introduction}\label{sec:intro}

In the doctoral register of the Faculty of Philosophy in Budapest, the doctorate of \enquote{margittai Neumann János} is entered under the academic year 1925/26. His dissertation, \emph{Az általános halmazelmélet axiomatikus felépítése} [The axiomatic construction of general set theory], bears the note \enquote{Elfogadták: Fejér Lipót és Suták József}: accepted by Lipót Fejér and József Suták.\footnote{\label{fn:register}ELTE University Archives (ELTE Egyetemi Levéltár), Doktori szigorlati jegyzőkönyv 1925/26, no.~26; reproduced in \citet{kaposi2023}. Kaposi does not mention Suták's name. The journal version of the dissertation is \citet{neumann1928}.} The first name needs no introduction. The second belongs to a Piarist priest, then sixty years old, who had written on the logical problems of set theory in 1911 and on the axiomatic construction of geometry in 1914, and who is scarcely more than a footnote in the history of Hungarian mathematics.

This article reconstructs the career of this man. József Suták (1865--1954) left no theorems bearing his name. Precisely for this reason, he allows us to study what the familiar narrative of the \enquote{Hungarian phenomenon} in mathematics, with its focus on the geniuses who later emigrated, easily overlooks \citep{pallo2004,frank2009}: how Budapest mathematics around 1900 received Bolyai and took up the foundational debates of the time; by what standards the faculty filled its chairs; and how a university teacher acted in the years of the numerus clausus. Three questions guide the account. First, what intellectual thread connects Suták's Bolyai translation of 1897 and his papers on set theory (1911) and on the axiomatics of geometry (1914) with his report on von Neumann's dissertation? Second, what does his appointment of 1911/12, when he was preferred to the Riesz brothers and to Haar, reveal about the faculty's priorities? Third, how should his role as chairman of the admissions committee after 1920 be assessed?

Research on this topic is limited. Besides the encyclopedia entries \citep{szinnyei1909,tudosnaptar} and the tribute in the annual report of his \emph{gymnasium} \citep{ertesito1912}, there are individual findings: on his teacher Gyula Vályi \citep{olahgal2014}, on the Eötvös Collegium \citep{garai2016} and, for the period after 1920, Géza Komoróczy's study of the faculty records on the numerus clausus \citep{komoroczy2015}. A coherent account is lacking.\footnote{This text is a thoroughly revised and considerably expanded version of an earlier English preprint (arXiv:2603.22196v1), which it corrects on several points on the basis of the sources cited below. A German version has been submitted to \emph{SieB -- Siegener Beiträge zur Geschichte und Philosophie der Mathematik}.}

The account rests on three bodies of sources. First, Suták's printed writings, including the offprints of 1911 and 1914. Second, records of the ELTE University Archives (faculty minutes, doctoral registers) and the archival material analysed by Komoróczy and by \citet{garai2016}. Third, a previously unknown part of Suták's personal papers: diplomas, appointment and salary records, correspondence, his will, and a typewritten obituary from the Piarist order. These items are in the author's private possession; they came to him through the Kosaras family, to which Suták's heir, Irma Kosaras née Kovács, belonged.\footnote{Cited below as \enquote{Suták Papers}. The obituary (typescript, 4~pp., after July 1954, author unknown) belongs to a genre with edifying intent and contains corrections in a second hand; it is used only where its statements are supported elsewhere or are explicitly marked as testimony.}

\section{Origins and education (1865--1891)}\label{sec:origins}

József Suták was born on 5~November 1865 in Szabadka (today Subotica, Serbia), in Bács-Bodrog County. The order's obituary records that he described himself as a Bunjevac, that is, as a member of the Catholic, South Slavic-speaking population of the Bácska, which had been distinguished by Hungarian national sentiment since 1848.\footnote{Obituary, p.~1 (Suták Papers): \enquote{Ő maga bunyevácnak mondta magát.}} After six years of gymnasium in Szabadka he entered the Piarist order on 27~August 1883; he was admitted by the Provincial Endre Palmár, a relative from Szabadka. The novitiate in Vác (1883/84) was followed by the last two years of gymnasium at the order's school in Kecskemét, where he passed the school-leaving examination on 5~June 1886, and by two years of theology in Nyitra (1886--1888).\footnote{\citet[5]{ertesito1912}; school-leaving certificate, obituary, p.~1 (Suták Papers).} He was ordained priest on 14~July 1889.

For his teacher training the order sent him to the University of Kolozsvár (Klausenburg, today Cluj-Napoca), where he studied mathematics and physics from 1888 to 1890 and passed the two teacher's examinations in May 1889 and May 1890. After a year at the order's gymnasium in Szeged (1890/91), during which he received his teaching diploma on 14~May 1891, he moved for the school year 1891/92 to the Piarist gymnasium in Budapest, where he taught until 1912. His dissertation on the curves of isoclinal normals on surfaces of the second order, supervised by Gyula Vályi, appeared in print in 1891 \citep{sutak1891,olahgal2009}. The doctoral diploma itself is dated 18~March 1892; it confers on him, already styled \emph{professor gymnasii ordinarius}, the degree of Doctor of Philosophy \emph{summa cum laude} in \emph{mathesis sublimior et elementaris}, and bears the signature of the rector, Antal Koch.\footnote{Doctoral diploma, Kolozsvár, 18~March 1892, reg.\ no.~811 (1891/2) (Suták Papers). This corrects \citet{szinnyei1909}, who dates the teacher's examination and the doctorate together to 1890.}

Vályi (1855--1913) remained a point of reference for Suták. When a chair of mathematics was advertised in Budapest in 1898 after the death of Gusztáv Kondor, Suták, as his \enquote{grateful pupil}, wrote to his teacher asking him to apply---evidence not only of his attachment to Vályi but also of the young gymnasium teacher's interest in appointments in the capital.\footnote{Letter from Suták to Vályi, Budapest, 8~February 1898, edited in \citet[4]{olahgal2014}: \enquote{Nagyságodnak kiváló tisztelője s hálás tanítványa}.} Kolozsvár in these years was not a uniform milieu of Bolyai veneration. Sámuel Brassai considered non-Euclidean geometry a chimera and believed as late as 1897 that he had proved the parallel postulate \citep[4]{olahgal2014}; Vályi is regarded as the \enquote{influential continuator} of Mór Réthy's work on Bolyai in Kolozsvár. However, he began his special lectures on Bolyai's \emph{Appendix}, about a third of which were devoted to its prehistory and which then commented on the text paragraph by paragraph, only in the second semester of 1891/92 \citep{szenassy1978}.\footnote{\citet{szenassy1978}: \enquote{Réthy munkájának a kolozsvári egyetemen nagyhatású folytatója volt Vályi Gyula (1855--1913), aki az 1891/92. tanév második félévétől kezdve többször tartott speciálkollégiumot Bolyai János Appendixéről.}} That Suták became acquainted with the \emph{Appendix} in Vályi's seminar therefore cannot be documented; an exchange in connection with his doctorate, which falls precisely in that academic year, is conceivable.

\section{Gymnasium teacher and Privatdozent (1891--1912)}\label{sec:teacher}

For two decades Suták was first and foremost a gymnasium teacher. He taught mathematics at all levels and, from 1905, physics in the upper classes; for decades he looked after the pupils' library, whose catalogue he published in 1901; and he wrote textbooks: a \emph{Számtan} (Arithmetic) for the lower classes (1898; six editions by 1908) and a two-volume \emph{Algebra} for classes IV to VIII (1901), for which he signed a contract with the publisher Wodianer in November 1900.\footnote{\citet{szinnyei1909}; publishing contract of 27~November 1900 (Suták Papers).} Former pupils recalled that they had learned from his lecture notes rather than from the prescribed textbook.\footnote{Letter from József Aurél Vécsey to György Balanyi, 27~July 1954, quoted in the obituary, p.~2 (Suták Papers).}

Alongside this, Suták built an academic career. His application for habilitation at the Faculty of Philosophy of the University of Budapest was considered at the faculty meeting of 20~February 1896. In the colloquium on 24~February a committee consisting of the dean Emil Thewrewk, Izidor Fröhlich and Ágost Scholtz examined him on the three fundamental theorems of Lie's theory of transformation groups, on minimal curves and isometric coordinate systems, and on the projective derivation of metric concepts and the role of imaginary elements. He gave his trial lecture on the singularities of algebraic curves. The habilitation was confirmed by ministerial decree no.~12.642 of 18~March 1896.\footnote{Habilitation records (Suták Papers). Radó Kövesligethy and Géza Bartoniek attended the trial lecture but were not members of the committee.} In the following 32 semesters Suták lectured as a Privatdozent for three to four hours a week, mainly on the theory of algebraic equations, Galois theory, invariants, algebraic curves, elimination and substitutions, each time before forty to ninety registered students. In 1898--1900 he served as a substitute lecturer. In the academic year 1910/11 he lectured on \enquote{set theory} to 68 students and on its applications to 52.\footnote{Certificate of the quaestor of 28~April 1911 on the lectures 1896--1911 (Suták Papers); \citet[7]{ertesito1912}.}

Suták was also involved early in teacher training. From 1895 he was a member of the national examination board for secondary-school teachers. In the same year Loránd Eötvös appointed him one of the first four subject tutors at the newly founded Eötvös Collegium, which, on the model of the École normale supérieure, was to train future gymnasium teachers; together with Géza Bartoniek, Sándor Mika, Frigyes Riedl and Móric Szilasi he laid down the Collegium's course of study. In 1898 he resigned from his tutorship, citing his workload \citep[pp.~50, 70, 74, 472\,f.]{garai2016}.\footnote{Invitation from Eötvös of 12~September 1895, with an annual fee of 600~florins (Suták Papers); ministerial appointment 50119/1895 of 17~September 1895 \citep[50, 70]{garai2016}. Garai stresses that \enquote{very few data} exist on Suták's activity at the Collegium \citep[74]{garai2016}.}

The scientific work of these years shows a mathematician of broad interests. Besides algebraic investigations in function theory (1894) and a new proof of a theorem of Eisenstein (1894), in 1895 he presented a \enquote{General theory of surfaces of the second class}, which developed their properties independently instead of deriving them, as was usual, from those of surfaces of the second order by the principle of duality \citep{sutak1895a}. For a monograph on the electromagnetic theory of light he received the Krisztina Lukács Prize of the Hungarian Academy of Sciences in 1895 \citep{sutak1895b}; the accompanying medal survives.\footnote{Medal by Joseph Daniel Böhm, obverse with the legend \enquote{Borúra derű}, reverse \enquote{M.\,T.\,A. Suták Józsefnek 1895} (Suták Papers).} There followed a 388-page theory of the differential and integral calculus (1900), a report on recent research on electricity (1904) and a theory of differential equations (1906) \citep{sutak1900,sutak1904,sutak1906}. The annual report of his gymnasium called the textbook of 1900 \enquote{gap-filling}, since no independent Hungarian work of this scope had appeared for twenty years, and that of 1906 the first to treat differential equations systematically \citep[8]{ertesito1912}. This is the judgement of his fellow members of the order, not of professional reviewers; but it shows what Suták's reputation rested on: the ability to present whole fields in an orderly way to Hungarian readers.

\section{Bolyai in 1897 and the geometric axioms of 1898}\label{sec:bolyai}

In 1896, the year of the millennium celebrations, Bolyai's \emph{Appendix} (1832) had long been available in Italian, French, English and German translations, but not in Hungarian. In 1897 two Hungarian translations appeared at once: that of Ignác Rados in the \emph{Mathematikai és Physikai Lapok}, and Suták's as a separate volume with preface, translation, commentary and a biography of Bolyai by the Temesvár architect Ferenc Schmidt \citep{sutak1897}.\footnote{Rados's translation: \emph{Mathematikai és Physikai Lapok} 6 (1897), from p.~145. Suták's edition was published by Schmidt himself (\enquote{Kiadja Schmidt Ferencz építész}) on commission with the university bookseller Frigyes Kilián (xxviii + 143~pp.); the Latin text is reproduced in facsimile. At the same time the Academy began its lavish edition of Farkas Bolyai's \emph{Tentamen} together with the \emph{Appendix} (1897--1904).} George Bruce Halsted, himself a translator of the \emph{Appendix}, reviewed Suták's edition in \emph{Science} in 1898, together with the first volume of the Academy's edition of the \emph{Tentamen} \citep{halsted1898}. Barna Szénássy later judged Suták's version \enquote{the better one as a translation}.\footnote{\citet{szenassy1978}: \enquote{Az egyik Rados Ignác, a másik -- mint fordítás a jobbik -- Suták József munkája.}} That the Academy preferred Suták's version to Rados's, or that there was a controversy between the two, is not documented. Nor do the gaps and errors that later Bolyai scholarship found in Schmidt's early biography concern Suták's translation \citep{acs2004}. The annual report of his gymnasium put Suták's merit more modestly and more aptly: before the Academy published the works of the Bolyais, Suták had \enquote{made accessible to everyone the work that set off the more recent geometrical research}.\footnote{\citet[8]{ertesito1912}: \enquote{Mielőtt a M.\,T. Akadémia e hiányt a Bolyaiak műveinek kiadásával pótolta volna, Suták mindenki számára hozzáférhetővé tette az újabb geometriai kutatásokat megindító munkát.}}

A year later Suták himself described his translation process. He had rendered Bolyai's symbol for parallelism by \enquote{párhuzamos} (parallel), for we are free \enquote{to name this symbol, in other words to translate it}; for Bolyai's limit surface and limit line (horosphere and horocycle) he introduced the names \enquote{Bolyai surface} and \enquote{Bolyai line}; and his commentary contained a proof of the consistency of Bolyai's geometry after Beltrami.\footnote{\citet[16--18]{sutak1898}: \enquote{nekünk azonban szabadságunkban áll ezt a szimbólumot elnevezni, más szóval lefordítani}; \enquote{Ennek a bebizonyítása megtalálható abban a »Bolyai-féle Geometriában« is, melyet az Appendixhez magyarázatként írtam.}} The edition of 1897 was thus more than a translation: it read the text of 1832 through the concepts and results of the period after 1868---in terms of the distinction drawn by \citet{grattan2004}, as \enquote{heritage} rather than \enquote{history}.

On 23~January 1898 Suták, as an ordinary member of the scientific and literary section of the Catholic Saint Stephen Society (Szent-István-Társulat), lectured on \enquote{Geometric axioms}; the lecture appeared in the \emph{Bölcseleti Folyóirat} and as an offprint in the Society's lecture series \citep{sutak1898}. He states his programme at the outset. He wished to show that \enquote{mathematical science, too, is of an empirical character}, that the inductive premises placed at its head are by no means as easy to obtain as had long been believed, and that precisely this difficulty had produced the three geometries of which, \enquote{owing to the unreliability of experience}, we cannot decide which actually exists.\footnote{\citet[4]{sutak1898}: \enquote{kimutassam, miszerint a mathematikai tudomány is tapasztalati jellegű [\dots] s hogy épen ezen a téren felmerülő nehézségek hozták létre azt a három, egymástól különböző geometriai rendszert, melyekről épen a tapasztalás megbizhatlansága folytán nem tudjuk eldönteni, hogy tulajdonképen melyik létezik.}} The lecture follows this programme in three steps. A history of the parallel problem from Euclid to Bolyai and Lobachevsky relies explicitly on the collection of documents edited by \citet{staeckelengel1895}; the consistency of Bolyai's geometry, Suták notes, was shown not by Bolyai himself but only by Beltrami, supplemented by Klein and by Gyula Kőnig's real model of 1872. A second part places the empirical principles on which Farkas Bolyai had based geometry in the \emph{Tentamen} (continuity and infinity of space, congruence, free mobility, monodromy) alongside the hypotheses of Riemann and Helmholtz; a third treats the projective metric of Cayley and Klein. The conclusion calls Euclidean (parabolic), Bolyai's (hyperbolic) and Riemann's (elliptic) geometry \enquote{perfectly equivalent} systems. Following Klein, Suták defines axioms as \enquote{requirements} by means of which we rise from the imprecision of intuition to unlimited precision; their consequences must be tested against experience, as in the exact natural sciences. To the question which of the three geometries exists, the last sentence replies: \enquote{We cannot decide.}\footnote{\citet[pp.~31\,f.]{sutak1898}: \enquote{E három egymással tökéletesen egyenlő értékű geometriai rendszer}; \enquote{Az axiómák oly követelmények, melyeknek segítségével a szemlélet pontatlanságáról [\dots] a határtalan pontosságra emelkedünk}; \enquote{hogy a három geometriai rendszer közül melyik létezik [\dots]? Nem tudjuk eldönteni.}}

The lecture is a well-informed and up-to-date survey, not original research; the history is derivative, and there are minor slips.\footnote{Saccheri's \emph{Euclides ab omni naevo vindicatus} is dated 1773 instead of 1733 \citep[8]{sutak1898}; Lobachevsky's names for the limit line and the limit surface (horocycle, horosphere) are interchanged \citep[17]{sutak1898}. Suták considers mistaken Lie's demonstration that Helmholtz's requirement of monodromy is dispensable in space \citep[pp.~29\,f.]{sutak1898}; cf.\ \citet{lie1893}.} The national emphasis is striking: Riemann and Helmholtz, he claims, had obscured the foundations with their coordinates and, without suspecting it, had woven the laurel wreath around the brow \enquote{of a Hungarian mathematician perhaps not even known to them}, namely Farkas Bolyai.\footnote{\citet[21]{sutak1898}: \enquote{nem is hitték, hogy evvel egy általuk talán nem is ismert magyar mathematikus homloka körül fonják a babérkoszorút.}} What is missing is remarkable, too. A priest of a religious order declares before a Catholic learned audience that the validity of Euclidean geometry cannot be decided empirically, and says not a word about the consequences for one's worldview. The annual report of 1912 turned this into the thesis that the geometries of Bolyai and Riemann are \enquote{equally true} with Euclid's;\footnote{\citet[8]{ertesito1912}: \enquote{E munka eredménye, hogy a Bolyai, valamint a Riemann-féle geometriák egyenlően igazak az Euklides-félével.}} Suták himself spoke more cautiously of equivalence and empirical undecidability. Philosophically, in 1898 he stands in the empiricist line of Gauss, Riemann and Helmholtz, not with Poincaré's conventionalism: one of the three geometries exists; we just cannot know which. By 1914 he would move far away from this position.

\section{Logic, sets, axioms: the papers of 1911 and 1914}\label{sec:papers}

\subsection{The logical problems of set theory (1911)}\label{sec:1911}

Suták's paper \emph{A halmazelméletben föllépő logikai problémák} [The logical problems arising in set theory] appeared in 1911 in the \emph{Közlemények} (Proceedings) of the Hungarian Philosophical Society \citep{sutak1911}. It is closely connected with the course on set theory that Suták gave at the university in 1910/11, and it is explicitly addressed to a philosophical audience: since the research on the antinomies concerned \enquote{the shaking of principles hitherto held to be infallible and of our confidence in our thinking}, he wished to present this field, so far known only to mathematicians, together with its \enquote{sources} \citep[2]{sutak1911}. In Budapest the debate was not new: Gyula (Julius) Kőnig had taken a position on the continuum problem at the Heidelberg congress of 1904 and in the \emph{Mathematische Annalen} in 1905 \citep{koenig1905}; his \emph{Neue Grundlagen der Logik, Arithmetik und Mengenlehre} followed in 1914 \citep{koenig1914}. The term \enquote{foundational crisis} (\emph{Grundlagenkrise}) goes back only to Weyl (1921) and is avoided here for 1911; to speak of a \enquote{debate on the antinomies} describes the situation more accurately \citep{weyl1921,ferreiros1999}.

The paper opens with a programme. Mathematics, Suták writes, is \enquote{a model built from elements constructed by reason according to the rules of logic}, by means of which we come to know the real world through comparison; for precisely this reason, \enquote{not the slightest deviation from the laws of logic} may be permitted in its construction. If a contradiction nevertheless appears, the part in question should not be torn down; rather, one must seek the \enquote{sources} from which the irregularity springs.\footnote{\citet[1]{sutak1911}: \enquote{A matematika az ész által konstruált elemekből a logika szabályai szerint fölépített oly modell, [\dots] misem természetesebb, mint az, hogy fölépítésében a logikai törvényektől a legcsekélyebb eltérés sem engedhető meg.}} By this standard Suták distances himself from two sides. Russell and Poincaré, whose works he cites up to the theory of types of 1908 and the Göttingen \emph{Sechs Vorträge} of 1910, wanted to steer thought, \enquote{by modifying the laws of the relations of thought}, in a direction in which the sources of the contradictions no longer occur; such a restoration of the disturbed equilibrium was, in his view, \enquote{not only not possible but also not justified}. Frege, Bernstein, Jourdain, Hessenberg, and Grelling and Nelson, on the other hand, had done \enquote{work that has so far been fruitless}, because none of them had traced the sources closely enough.\footnote{\citet[pp.~1\,f.]{sutak1911}: \enquote{gondolkodásunk megzavart egyensúlyának ily módon való helyreállítása nemcsak hogy nem lehetséges, de nem is jogosult.} Cf.\ \citet{russell1908,poincare1906,poincare1910}.}

Suták's own starting point is a naive, content-based notion of set: a set is the totality of the elements that satisfy certain criteria. From this he derives a remarkably strong existence principle. If the criteria contain no contradiction, \enquote{the set exists, whether or not we can decide of the individual things whether they are elements of the set}. His examples are the algebraic numbers, of which one does not know whether $2^{\sqrt2}$ belongs to them, and the number $\pi$, which did not belong to them even before Lindemann's proof of 1882 \citep[pp.~3\,f.]{sutak1911}. The decidability of membership is thus irrelevant to existence, a view that may be called realist. Consistently, Suták also admits the set of all things that \enquote{contain themselves as an element} \citep[4]{sutak1911}. His formulation of the basic logical principle is noteworthy: as the \enquote{principle of contradiction} he states \enquote{every thing is either $a$ or not $a$} \citep[3]{sutak1911}, that is, what is today called the law of excluded middle. The conflation of the two principles has consequences for his argument, for it allows him to treat every consistently defined property as a division of the whole universe into two sets.

On this basis Suták distinguishes two sources of the antinomies. The antinomies of Russell and of Burali-Forti stem, in his view, from definitions that themselves contain a contradiction or contradict a theorem of set theory \citep[24]{sutak1911}. Richard's antinomy, by contrast, rests on a confusion. Let $A$ denote the set of concepts already designated in language and $G$ the set of decimal fractions definable by finite combinations of these concepts. As soon as $G$, too, is given a name and added to $A$, an enlarged set of concepts $A'$ arises, with a set $G'$ of definable decimal fractions that contains $G$ but may contain new elements. The antinomy springs from the \enquote{erroneous supposition} that $G$ must also contain the elements of $G'$; its second form springs from the \enquote{tacit exchange} of $G$ for the totality $H=(G,G',\dots)$ of the decimal fractions definable by concepts already constructed \emph{and still constructible} \citep[25]{sutak1911}. The paper concludes:
\begin{quote}
The set-theoretical contradictions thus arise partly from the tacit exchange of sets of different content, partly from definitions containing a contradiction, but not from a limited validity of the logical principle of contradiction.\footnote{\citet[26]{sutak1911}: \enquote{A halmazelméleti ellentmondások tehát részint különböző tartalmú halmazok szubkonszciens módon való felcseréléséből, részint pedig ellentmondást tartalmazó definiciókból-, nem pedig a logika ellentmondási elvének korlátos érvényességéből erednek.}}
\end{quote}

In substance, the analysis of Richard's antinomy is close to that of Richard himself and of Poincaré, who likewise located the error in the fact that the totality of definable numbers is not yet closed while the definition is being given \citep{richard1905,poincare1906}. In his programme, however, Suták is their opponent: Poincaré drew consequences from this for the admissibility of impredicative definitions and thus for logic itself, whereas Suták holds fast to classical logic without restriction and to the content-based notion of set. Assigning him to the later \enquote{schools} of intuitionism or formalism would be anachronistic for 1911.

The paper contains further sections that go beyond the antinomies and are revealing for its interpretation. Suták treats the well-ordering of the continuum as an open question; as the only \enquote{significant result} he cites Kőnig's proof that the cardinal number of the continuum cannot be among the alephs with index not less than $\omega$ \citep[21]{sutak1911}---a rendering that goes beyond what Kőnig actually proved (that the continuum is not $\aleph_\omega$) \citep{koenig1905,moore1982}. He calls the trichotomy of infinite cardinal numbers \enquote{still unsolved} \citep[22]{sutak1911}; Zermelo's well-ordering proof of 1904, from which it follows, is not mentioned in this context, although Suták cites Zermelo elsewhere \citep{zermelo1904}. Finally, he interprets complete induction neither as a definition nor as a logical principle but as a method of proof showing that a theorem is independent of the order of a discrete set; for \enquote{singular} places such as $\aleph_0$ in the sequence $1,2,3,\dots,\aleph_0$ he introduces a \enquote{hyperinduction}, and from there rejects Poincaré's view that the consistency of the definitions placed at the head of a science can be secured only by complete induction \citep[22--24]{sutak1911}.

The order's report of 1912 summarized the result approvingly, emphasizing that Suták had put forward \enquote{a new procedure of original conception} for removing individual antinomies \citep[8]{ertesito1912}. No other contemporary reactions are known. The paper is therefore of interest less as a contribution to the international debate than as a document of how a well-informed Budapest mathematician dealt with the antinomies in 1911: with precise knowledge of the literature, with a strong, content-based and realist notion of existence, and with the conviction that it was not logic but our concept formation that needed correcting. Theological arguments do not occur in the text.

\subsection{The logical construction of geometry (1914)}\label{sec:1914}

Three years later Suták, now a full professor entrusted with teaching analytic and synthetic geometry, published in the same journal the paper \emph{A geometria logikai épületének felépítése} [The construction of the logical edifice of geometry] \citep{sutak1914}. It is his only known work in which he himself erects an axiomatic system.

The preface engages with the contemporary axiomatics of projective geometry. Suták mentions Fano's postulate (if $C$ and $D$ are harmonic with respect to the distinct points $A$ and $B$, then $C$ and $D$ are also distinct), its form in Veblen and Young (the diagonal points of a complete quadrangle are not collinear), Pieri's definition of the segment by means of harmonic separation (with reference to Russell and Whitehead), and Veblen and Young's postulate of projectivity \citep[pp.~1\,f.]{sutak1914}.\footnote{Cf.\ \citet{fano1892,pieri1898,russell1903,whitehead1906,veblenyoung1910}; Suták dates Fano's paper to 1891. Hilbert is not mentioned \citep{hilbert1899}.} Such postulates, he says, are \enquote{logically unobjectionable} but presuppose a thorough knowledge of the finished edifice and do not meet \enquote{the most elementary requirements of simplicity}. He himself uses no postulates \enquote{that belong to the higher regions of mathematics}; the only undefined primitive notions are \enquote{points} and \enquote{bipartition} (\emph{kettéválasztás}), of which one knows only that certain sets are divided into two subsets by two of their elements \citep[2]{sutak1914}. He explicitly refrains from proving the \enquote{existence principle}, that is, the consistency of the system, \enquote{since this is not possible on a purely logical basis}. Neither recourse to sets of empirical elements nor recourse to sets of numbers provides such a proof: since sets of points and sets of numbers are isomorphic, the existence principle presupposed for the set of numbers by no means justifies that for the set of points \citep[pp.~2\,f.]{sutak1914}. In a footnote he describes the method of models: the compatibility of the postulates is shown by an interpretation of the undefined symbols that satisfies all postulates, their independence by interpretations that satisfy all but one \citep[3]{sutak1914}.

The introduction develops a brief epistemology to go with this. The primitive notions of geometry are prompted by experience, but geometric figures have no material existence; geometry is therefore to be regarded \enquote{as a creation of reason}, differing from other sciences only in needing fewer stimuli from experience. \enquote{Logicism}, which grew out of set theory, constructs it purely logically, yet cannot deny its origin: it always chooses its contents with a view to the results of experience \citep[3--5]{sutak1914}. Accordingly, Suták distinguishes \emph{pure} geometry, the logistically constructed edifice that has nothing to do with the intuition of points and space, from \emph{applied} geometry, which is built from empirical material and presupposes intuition; drawings serve only for orientation \citep[5]{sutak1914}. Logically identical edifices that differ only in terminology (dual ones) lose their identity as soon as their bearers of content are replaced by empirical elements \citep[24]{sutak1914}.

The construction itself is an order-based axiomatics of projective geometry. Five axioms of the line, governing incidence and bipartition and already implying that every line contains at least countably infinitely many points, are followed by the notions of the separation of point pairs, of cyclic order and of the segment, a Dedekind continuity axiom with the existence of limit points, and a study of ordered correspondences, whose double-point theorems Suták explicitly compares with the motion of two points on a circle \citep[6--16]{sutak1914}. With an axiom of the plane and an axiom of space he obtains plane and space, plane duality, duality in the bundle and spatial duality, and Desargues's theorem, which he proves from the axioms of space; via the complete quadrangle he arrives at harmonic points. The goal is the theorem that harmonic pairs separate each other and, from it, that the three diagonal points of a complete quadrangle are not collinear \citep[pp.~31\,f.]{sutak1914}. Fano's postulate thus becomes a theorem, in both of the forms mentioned in the preface. Methodologically, the text stands in the tradition of Pieri and of Enriques, to whose \emph{Lezioni} Suták refers several times \citep{enriques1898,marchisotto2007}. An announced continuation, which was to prove \enquote{the intuitive theorem of the Pythagoreans that everything is number}, did not appear in this form \citep[2, 28]{sutak1914}.

For an assessment of Suták, the paper matters in three respects. First, it shows the step from survey (1898) to construction of his own, and with it a changed epistemology: in 1898 Suták had called mathematics \enquote{empirical} and left the question of the actual geometry to experience; in 1914 pure geometry is a creation of reason, to which experience gives only the impulses. Second, it takes up the question of existence from 1911, but with a shifted emphasis: in 1911 a consistent criterion sufficed for the existence of a set; in 1914 Suták states that the consistency of an axiomatic system cannot be proved purely logically, and he distrusts arithmetization as a way out. Third, the text shows that Suták was familiar with the axiomatic method and its use of models. This makes it understandable why the faculty appointed him, alongside Fejér, as referee for von Neumann's axiomatic set theory in 1926 (Section~\ref{sec:professor}).

\section{The appointment of 1911/12}\label{sec:appointment}

At its last meeting of the academic year, on 1~June 1911, the Faculty of Philosophy had to decide on filling the third chair of mathematics. Six applications had been received: from Frigyes Riesz and Pál Dienes, both gymnasium teachers, from the Privatdozenten Suták and Mihály Demeczky, from the Göttingen Privatdozent Alfréd Haar, and from Marcel Riesz. The committee---Izidor Fröhlich, Loránd Eötvös, Manó Beke as rapporteur, Radó Kövesligethy and Jenő Klupathy---could not agree. The majority (Fröhlich, Eötvös, Beke) proposed Frigyes Riesz as extraordinary professor in first place, Suták as full professor in second place and Dienes as extraordinary professor in third place; the minority (Klupathy, Kövesligethy) placed Suták as full professor first, Haar second and Dienes third.\footnote{ELTE University Archives, minutes of the Faculty of Philosophy 1910/11, seventh (closing) meeting of 1~June 1911, dean's no.~1617/1910--11, item~10, no.~1560, fols.~242--248, here fols.~242\,f. The minutes refer to the \enquote{coming academic year 1911/12}. A detailed report by Beke on each candidate was read out but is not in the file.}

In the discussion the art historian Gyula Pasteiner took the floor. The committee, he said, had treated the university's own Privatdozenten \enquote{in a very step-motherly fashion}; it should have acknowledged that Demeczky was a good teacher who always fulfilled his duties punctually and had worked irreproachably for eighteen years. What moved him was \enquote{the principle of the solidarity of university teachers}. The Anglicist Arthur Yolland agreed: Pasteiner had rightly stressed that teacher training was \enquote{an important factor} of the faculty, and Demeczky was \enquote{most excellent} as a teacher; he moved the ranking Suták, Demeczky (both as full professors), Dienes (extraordinary).\footnote{Ibid., fols.~243\,f.: \enquote{a bizottság egyetemünk magántanáraival nagyon mostohán bánt el}; \enquote{az egyetemi tanári szolidaritás elve}; \enquote{a karnak fontos tényezője a tanárképzés}.} In a secret ballot Suták received 24 of 41 votes for first place, Riesz 11, Demeczky 4. For second and third place no ballot produced an absolute majority; in the run-offs Demeczky and Riesz stood at 16 to 13 with eleven blank ballots, and Demeczky and Dienes at 18 to 18. The faculty resolved unanimously to propose Suták as full professor. Pasteiner and Yolland had a minority opinion in favour of Demeczky entered in the minutes.\footnote{Ibid., fols.~246--248.} The faculty thus proposed a single candidate.

The ministry's decision differed from what the proposal would have led one to expect. By royal resolution of 21~March 1912, \emph{Suták and Demeczky} were both appointed full professors. Demeczky, a Privatdozent since 1893, was at that time director of the Franz Joseph Educational Institute in Budapest and seconded to the Ministry of Religion and Public Education.\footnote{Appointment decree no.~42119 of 21~March 1912 and salary decree no.~60196/1912 (Suták Papers).} Riesz, Haar and Dienes came away empty-handed; Riesz obtained a professorship in Kolozsvár in 1911, Haar in 1912. Suták took the oath of office before the rector Izidor Fröhlich on 20~April 1912 and received a salary of 6000~crowns from 1~May 1912; the order's permission to apply and to work outside the order dated back to 3~May 1911. A decree of June 1913 assigned him analytic and synthetic geometry, and Demeczky number theory and algebra, each with at least five hours a week.\footnote{Record of the oath, salary order, the order's permission of 3~May 1911 and decree no.~79953 of 10~June 1913 (Suták Papers).}

In retrospect it may seem strange that the faculty placed Suták ahead of Riesz and Haar. The minutes, however, allow the reasons to be stated more precisely than talk of a wrong decision does. In Suták's favour was the fact that the committee majority, too, proposed him for a full professorship, whereas it proposed Riesz only for an extraordinary one; also in his favour were fifteen years of teaching as a Privatdozent and, in the logic of Pasteiner and Yolland, the faculty's obligation towards its own Privatdozenten and teacher training as one of its tasks. That the ministry appointed Demeczky, of all people, the candidate of the minority opinion, shows how strong these considerations were.

There is no mention in the minutes of a \enquote{philosophical wonder} that the faculty might have valued in Suták. This formula comes from another source. On his departure, the annual report of his gymnasium devoted a tribute to him in 1912 that reveals the order's self-image \citep{ertesito1912}. It asks whether in Suták \enquote{the scholar or the educator is the greater} and answers in favour of the \emph{praeceptor Scholarum Piarum}. It quotes Aristotle's remark in the \emph{Nicomachean Ethics} (VI~9) that a child may become a mathematician but not wise, and Descartes's \emph{omnia apud me mathematice fiunt}. It calls his inexhaustible intellectual energy that of a \enquote{radium man} (\emph{rádiumember}) and praises his gift of making others think, \emph{faire penser les autres}. And it characterizes him as a thinker by the \enquote{wondering immersion of the mind in problems}, \enquote{roughly what the Greeks called} θαυμασμὸς φιλοσόφος, combined with a \enquote{certain logical tenacity}.\footnote{\citet[5--9]{ertesito1912}, quotations on pp.~5, 7, 6 and 9: \enquote{Nehéz volna eldöntenünk, vajjon Sutákban a tudós vagy a pedagógus-e a nagyobb?}; \enquote{a szellemnek a problémákba való csodálkozó elmerülése jellemzi; körülbelül az, amit a görögök θαυμασμὸς φιλοσόφος-nak neveztek.} The Greek phrase is printed thus; presumably Plato's \enquote{μάλα γὰρ φιλοσόφου τοῦτο τὸ πάθος, τὸ θαυμάζειν} (\emph{Theaetetus} 155d) is meant.} This characterization is a testimony of the order's perspective, not a justification given by the faculty. As such, however, it is revealing, because it describes the mathematician, precisely in his effect as a teacher and in his contemplative seriousness, as the embodiment of the Piarist ideal.

\section{Professor in Budapest (1912--1936)}\label{sec:professor}

As a full professor Suták taught mainly analytic and synthetic geometry. The recollections of former students from the order, collected in the obituary, portray a teacher whose effect lay in logical rigour and an almost theatrical enthusiasm. János Pogány wrote that Suták had presented analytic geometry \enquote{so synthetically and deductively and so logically constructed in all its parts} that he had never found anything like it since in any book. In moments of enthusiasm Suták would personify the concepts and, striking the blackboard with his fist, call out: \enquote{You, point, go to infinity, and you, line, turn about the point!} The students positively sought out his lectures on the infinite, because he also touched on theology; after a lecture on point, line and plane, according to Ödön Hesz, the listeners concluded: \enquote{He spoke about the Trinity.} He rejected the teaching of determinants in the gymnasium as mere formalism that stifled independent thinking.\footnote{Obituary, pp.~4\,f. (Suták Papers): \enquote{Az analitikus geometriának olyan szintetikus és deduktiv és minden részében olyan logikusan felépitett módszerét adta, amelyet azóta sem találtam meg egyetlen e tárgykörrel foglalkozó matematikai könyvben sem}; \enquote{Te pont, te menj a végtelenbe, te egyenes pedig forogj a pont körül!}; \enquote{A Szentháromságról beszélt.}} Such testimonies, especially in an obituary, must be read with caution. But they show that the theological dimension one looks for in vain in the papers of 1911 and 1914 had its place in the lecture hall.

A testimony from outside the order comes from the writer Géza Ottlik, who studied mathematics and physics in Budapest in the early 1930s. In a recollection published in 1986 he names Suták, together with Fejér and Ortvay, among the teachers who compensated him for having given up the study of literature. Suták, he writes, was \enquote{an old secular priest} who lectured on projective geometry \enquote{with a weakening voice and unbroken intellectual freshness}; in a threadbare cassock he walked in from Nagytétény, always bringing his lunch in a little bundle tied in a kitchen towel, and he never failed to stress again and again \enquote{that mathematics need not be learned, because it is within us}.\footnote{\citet{ottlik1986}: \enquote{Suták öreg világi pap volt, projektív geometriát adott elő, gyengülő hangerővel és töretlen szellemi frissességgel. Kopott reverendájában gyalog járt be Nagytétényből, kis konyharuha-batyuban hozta magával mindig az ebédjét. Nem felejtette el újra meg újra hangsúlyozni, hogy a matematikát nem kell tanulni, mert az bennünk van.} The text is a recollection given literary form half a century later; details such as the walk from the suburb are not to be taken literally.}

At the same time Suták took on the usual academic offices. In the academic year 1913/14 he was notary of the faculty; from 1915 to 1920 he examined mathematics on the examination board for secondary-school teachers, and from 1922 to 1926 he was a member of the examination board for teachers at teacher-training colleges. With effect from 31~December 1918 he was placed in the fifth salary grade; the letter of appointment of February 1919 bears the signature of the prime minister Dénes Berinkey.\footnote{Appointment and salary records 1913--1926 (Suták Papers). The form of address \emph{Méltóságos} in later letters is connected with this rank; it denotes a rank in the higher civil service and does not correspond to \enquote{Excellency}.} In his scholarly work he now appeared mainly in the Saint Stephen Academy (Szent István Akadémia), founded in 1916. Between 1916 and 1928 ten of his lectures were published in its \emph{Értesítő} (Bulletin), some in full, some as a report or an excerpt.\footnote{\citet{sutak1916,sutak1917,sutak1920,sutak1921,sutak1922,sutak1923,sutak1924,sutak1925,sutak1926,sutak1928b}. The lectures of 1916 and 1925 also appeared with the Catholic publisher Stephaneum \citep{sutak1918,sutak1928}.} Most of them belong to differential geometry and use vectorial methods: from an \enquote{absolute theory of curves} (1916) and the \enquote{fundamental theorems of the theory of surfaces} (1919), via ruled surfaces and associated curves on surfaces, to the hypersurfaces of $(n+1)$-dimensional and the $n$-dimensional surfaces of $(n+m)$-dimensional Euclidean space (1924/25). In addition, there are a theorem on perspective tetrahedra following on from Desargues's theorem (1922) and an algebraic paper on collineations and congruences of higher degree (1923). For his \enquote{research on the theory of surfaces} from the years 1917 to 1925 he received the Academy's Fraknói Grand Prize.\footnote{\citet[12]{viszota1940}: \enquote{Suták József: A felületelméletre vonatkozó kutatások. Az 1917--1925. évkörből, a mennyiségtani osztályból.} The prize was awarded in rotation among the sections for periods of eight years each.} His last lecture (1928) was devoted to the claim, going back to Schläfli, that every $n$-dimensional Riemannian space is a subspace of a Euclidean space of suitable dimension; Suták wanted to give it a mathematical foundation because he had found no acceptable one. In three pages, using Riquier's existence theory for passive systems of partial differential equations and Ricci's notion of the class of a Riemannian space, he sketches an argument that a Euclidean space of dimension $\binom{n+1}{2}$ always suffices.\footnote{\citet[49--51]{sutak1928b}: \enquote{Mivel a jelezett tételnek elfogadható matematikai indokolását nem találtam}; \enquote{a legkedvezőtlenebb esetben $r=\binom{n}{2}$, mert akkor $N'=n+r=\binom{n+1}{2}=N$.} Suták refers to Méray, Riquier (1893), Ricci (1884) and to his own work on differential equations of 1906 and 1907 \citep{sutak1906}.} This is exactly the bound for which Maurice Janet (1926) and Élie Cartan (1927) had just proved local embeddability; Suták does not mention them \citep{schlaefli1873,janet1926,cartan1927}. The example shows the possibilities and the limits of these late works: they took up current questions with appropriate means, but they appeared in Hungarian in a journal hardly noticed outside Hungary, and they did not always keep pace with the international literature. Outside the Academy, the enlarged second edition of his textbook of analysis appeared in 1922.

Suták belonged to the Academy as an ordinary member from the beginning, and from 1926 as an honorary member. After the death of the previous chairman, the geologist Ferenc Schafarzik, the section for mathematics and natural sciences elected him its chairman in December 1927 \enquote{for the remaining period}; he was re-elected in February 1929. In February 1938 he resigned \enquote{in view of his advanced age} and proposed the astronomer József Wodetzky as his successor. He remained an honorary member until the Academy was dissolved in 1951.\footnote{The volume for 1925 lists Suták throughout as an ordinary member \citep{szia1925}; the membership list of 1926 names him as the only honorary member of the fourth section, whose chairman was then Schafarzik \citep[99]{szia1926}. Election in December 1927 \citep[88]{szia1927}: \enquote{az osztály a hátralevő időre Suták József t.~tagot az osztály elnökévé választotta}; re-election at the general meeting of 22~February 1929 \citep[55]{szia1929}; resignation, announced at the closed meeting of 18~February 1938 \citep[121]{szia1938}: \enquote{Suták József tiszteleti tag, osztályelnök magas korára való tekintettel arra kéri az osztály tagjait, hogy [\dots] személyét kapcsolják ki a választásból.} The dates 1929--1938 in the obituary, p.~5 (Suták Papers), correspond to the term after his re-election; the \emph{Magyar Katolikus Lexikon} (art.~\enquote{Szent István Akadémia}) gives 1928--1937. The obituary wrongly has his honorary membership follow his chairmanship. On the celebration of his seventieth birthday at the meeting of 29~November 1935, which he himself chaired, see \citet[99]{horusitzky1939}.} The priest-scholar thus moved in two public spheres: the state university and the Catholic scholarly organization.

The episode mentioned at the outset also falls in these years. John (János) von Neumann enrolled at the Faculty of Philosophy on 14~September 1921, with mathematics as his major subject and experimental physics and chemistry as minor subjects; his record book lists Suták's geometry lectures alongside Fejér's analysis and Gusztáv Rados's algebra \citep[3]{szabo2003}. In parallel with his chemistry studies in Berlin and Zurich, he completed his Budapest studies on 11~July 1925. His dissertation on the axiomatic construction of general set theory, whose journal version appeared in 1928 \citep{neumann1928}, was accepted by Fejér and Suták. In the doctoral examination (\emph{rigorosum}) Fejér examined in mathematics, Károly Tangl in experimental physics and a third professor in chemistry; all subjects were graded \emph{summa cum laude}. The conferral of the doctorate took place on 13~March 1926.\footnote{Doctoral register (see n.~\ref{fn:register}); \citet{kaposi2023}. The date of the rigorosum reads as 10~March 1926; \citet[3]{szabo2003} gives 12~March. The signature of the chemistry examiner (probably Gusztáv Buchböck) is not fully legible.} Suták was thus not an examiner in the rigorosum, but he was one of the two referees of the dissertation. After his papers of 1911 and 1914 he was well qualified for this: the critic of the antinomies, who in 1911 had wanted to save set theory through careful concept formation and in 1914 had pointed out the limits of a purely logical proof of consistency, assessed in 1926 an axiomatic set theory that avoided the antinomies precisely by distinguishing between totalities that can themselves be elements and those that cannot---the germ of the later distinction between sets and classes.

In 1936, at the age of 71, Suták retired. According to the obituary, he regretted having to leave the lectern; had they let him stay there, he joked, he would have proved that not all men must die.

\section{Under the numerus clausus (1919--1936)}\label{sec:nc}

The years after 1919 confronted the Budapest faculty with political decisions. After the fall of the Hungarian Soviet Republic, students were screened by vetting committees (\emph{igazoló bizottságok}) for their conduct during the revolutionary months \citep{ladanyi1979}; with Act XXV of 1920, Hungary became the first European country of the postwar period to introduce an admissions quota, which in its application targeted Jewish students \citep{kovacs2012}. At the Faculty of Philosophy a radical wing around the zoologist Lajos Méhely and the art historian Antal Hekler pressed for strict application \citep[56--67]{komoroczy2015}. Suták was involved in the implementation in two capacities: as rapporteur of the student vetting committee and as chairman of the admissions committee.

Suták himself was not affected by the interventions of the Soviet government. When the People's Commissariat of Education withdrew the teaching licence of 25 professors of the faculty on 26~March 1919, among them Ballagi, Pauler, Yolland and Suták's co-appointee Demeczky, followed a few days later by Fináczy, Hekler and Tuzson, he seems to have remained unaffected. On the contrary, he appears on the list of 24 professors whose lectures continued \enquote{unchanged} in the spring semester of 1919 \citep[pp.~200\,f., nn.~30 and 37]{turan2024}.\footnote{Communicated on 10~April 1919 by Tódor (Theodore von) Kármán as deputy head of section in the People's Commissariat (EÜK 354/1919, after \citealp[p.~201, n.~37]{turan2024}). The list also includes Beke, Fejér, Kövesligethy and Méhely; Beöthy, Fröhlich and Mauritz, who had initially been suspended, had meanwhile been readmitted. Suták was not a member of the faculty's reform committee set up in February 1919 \citep[p.~199, n.~28]{turan2024}.} He did not, however, sign the declaration by which the teachers from the religious orders at the nationalized schools were to renounce the clerical state and enter the service of the Soviet Republic.\footnote{Name index (A magyar piaristák 1918/1919-ben), in \citet[p.~435]{mpm2013}: \enquote{nem írta alá a nyilatkozatot (pedig a Tanácsköztársaság területén levő rendház tagja volt)}. Wording of the declaration: \enquote{Kijelentem, hogy az egyházi pályát elhagyom és a Tanácsköztársaság szolgálatába lépek}; cf.\ \citet[pp.~22\,f.]{szakal2013}; according to Szakál, about two thirds of the Piarist teachers signed.} Suták thus continued to teach under the Soviet government without joining it. Unlike Ballagi, whose harshness Komoróczy explains partly by the memory of his own suspension \citep[57]{komoroczy2015}, Suták's role after 1920 cannot be interpreted as a response to slights he had suffered.

An example of the first capacity is the faculty meeting of 25~October 1920, at which Suták presented the report of the vetting committee. It recommended four students for admission, among them Árpád Illyefalvi Vitéz, against whom objections had been raised because he had been a member of the Galileo Circle. Since he was known, however, as a gifted and industrious young man who had been unsettled only temporarily by his severe material circumstances, the committee wished \enquote{to take the standpoint of forgiveness}. Hekler objected: for a faculty that trained teachers, character was \enquote{far more important even than talent}; he moved a resolution in principle no longer to admit members of the Galileo Circle. Suták defended Illyefalvi, citing in his favour that he had dared to stand up to the terror during the Soviet Republic. The faculty rejected his admission by 12 votes to 10.\footnote{ELTE University Archives, minutes of the Faculty of Philosophy 1920/21, meeting of 25~October 1920, item~31 (\enquote{A hallgatókat igazoló bizottság jelentése}), pp.~18\,f., fols.~52\,f.: \enquote{a bizottság vele szemben a megbocsátás álláspontjára kíván helyezkedni}; \enquote{a tehetségnél még sokkalta fontosabb a jellem}; \enquote{a kommunizmus idején föl mert lépni a terrorral szemben}.}

In the second capacity a different picture emerges. As chairman of the admissions committee, Suták took part in the decisions by which the faculty enforced the quota against the students of the National Rabbinical Seminary, who needed a doctorate from the Faculty of Philosophy for rabbinical ordination. Komoróczy has reconstructed the events from the faculty records \citep[62--67]{komoroczy2015}:
\begin{itemize}
\item On 8~March 1921 the committee under Suták's chairmanship found, with regard to the seminary's petition, which was supported by the ministry, that \enquote{there is no necessity for special measures}; the faculty concurred and struck off three seminarians who had already enrolled \citep[63]{komoroczy2015}.
\item On 9~September 1921 Suták presented the committee's report at an extraordinary faculty meeting: besides the Christian first-year students, 24 Jewish students were to be admitted, corresponding to 6~per cent of the number of 400 approved by the ministry, eight each from the Rabbinical Seminary, Jewish men and Jewish women; students wishing to transfer from other universities were to be vetted again. This was somewhat more than the 5~per cent that Méhely had demanded earlier \citep[64]{komoroczy2015}. At the same meeting, after contributions from Petz, Suták and the dean, the faculty declared that it did not wish \enquote{to deal with the matter on its merits} when the seminary asked that its students not be counted within the 6~per cent.\footnote{ELTE University Archives, minutes of the Faculty of Philosophy 1921/22, first extraordinary meeting of 9~September 1921, items~6 and 7, pp.~6\,f.: \enquote{A bizottság a keresztény hallgatókon kívül a minisztérium által engedélyezett 400-as létszám 6\%-ának megfelelőleg 24 zsidó hallgatónak, köztük 8 rabbiképző intézeti növendéknek, 8 zsidó férfihallgatónak és 8 nőhallgatónak felvételét javasolja}; \enquote{a Kar kimondja, hogy az üggyel nem kíván érdemlegesen foglalkozni}. The minutes do not record what Suták said in the discussion.}
\item On 6~April 1922 Suták presented a definition of the \enquote{nationality} referred to in the law: for \enquote{races or nationalities} whose character was closely bound up with religion, an extract from the birth register sufficed to establish it; only in the rare case of a transfer from another nationality should the declaration on conversion and descent be taken into account. Komoróczy comments sarcastically that it \enquote{took a mathematician's mind} to put the simple truth, \enquote{a Jew is a Jew}, into a complicated sentence \citep[63]{komoroczy2015}.
\item On 22~October 1925, in a draft resolution in his own hand, Suták declared the seminary's renewed petition acceptable, with the restriction that the seminarians might take no examinations at the faculty other than the doctorate \citep[65]{komoroczy2015}.
\item When the Minister of Religion and Public Education, Kuno Klebelsberg, decreed in November 1925 that the seminarians were not to be counted within the quota---against the background of the proceedings before the League of Nations---Méhely, Melich, Pauler, Suták and Tuzson protested on 26~November, invoking university autonomy, with \enquote{grave reservations} \citep[65]{komoroczy2015}.
\item In January 1926 Suták took part in the meeting at the ministry at which State Secretary Tóth asked the faculty to be accommodating, referring to the negotiations in Geneva; the restricted range of subjects for the seminarians had been negotiated beforehand in Suták's departmental rooms \citep[67]{komoroczy2015}.
\end{itemize}
John von Neumann, who enrolled as \enquote{izraelita} five days after the meeting of 9~September 1921 (Section~\ref{sec:professor}), was probably among the eight Jewish men who could be admitted in the autumn of 1921. The minutes do not reveal the criteria by which the committee chose among the applicants; besides \enquote{national reliability}, the law required selection according to intellectual ability \citep{kovacs2012}, and by this standard a school-leaver regarded as the best mathematics pupil in the country could hardly be turned away \citep[3]{szabo2003}. Komoróczy does not name Suták among the professors who, he assumes, voted in the minority against antisemitic motions \citep[80]{komoroczy2015}. The later doctoral examinations of the rabbinical candidates were conducted by other professors \citep[pp.~71\,f., 80--83]{komoroczy2015}; mathematics was not among the subjects permitted to them.

The Suták Papers add to this picture two letters showing the outside pressure to which the committee was exposed: an anonymous letter of 5~October 1921, directed, with an eye to the admissions committee, against the student Zsuzsa Back, and a recommendation from 1935 for an applicant named Bonta, which expressly mentions his Jewish origin---seven years after the amendment of 1928 had removed the explicit reference to \enquote{race and nationality} from the law.\footnote{Suták Papers. The outcome of both cases is unknown. The sender of the second letter is Antal Tass, director of the observatory and, like Suták, a member of the fourth section of the Szent István Akadémia \citep[99]{szia1926}.} Against this stand three episodes that the obituary relates from the recollections of students. In 1919, in defence of academic freedom of teaching, Suták is said to have confronted soldiers with fixed bayonets who forced their way into his lecture. During a student strike he went on lecturing to the few who were present, because there was no order from the rector. And when, later, antisemitic demonstrators burst into the lecture halls, demanding that teaching be stopped and the Jewish students removed, he is said to have turned them out with the words: \enquote{Get out of here, there is no Jew here!} The author of the obituary adds that Suták could hardly have ascertained at that moment whether there really was none in the room, and probably did not try to.\footnote{Obituary, p.~4 (Suták Papers), on the testimony of László Helyes: \enquote{Takarodjanak innen, itt nincs zsidó!} The episode is undated.}

These testimonies are heterogeneous and cannot be set off against each other. The faculty records show Suták as an agent who both implemented and shaped the quota, who made the category of Jewish \enquote{nationality} administratively workable and who in 1925 sided with the spokesmen of the radical wing against the ministry. The obituary shows him as someone who defended his lecture hall against interference of every kind, including antisemitic disrupters, while adopting their category in his defence. Common to both pictures is a standard that was institutional rather than political: order in teaching, the autonomy of the faculty, the validity of regulations. In individual cases this standard could protect, as in the case of Illyefalvi or in the lecture hall; but at the same time it made Suták a reliable executor of a discriminatory order.

\section{Later life and papers (1912--1954)}\label{sec:later}

On his appointment in 1912 Suták had had to leave the old house of the order on the Danube embankment, because his room was needed for his successor; even in old age he is said to have complained about it. He rented a flat with the Kovács family in what is today Bartók Béla út and took care of the two children, József and Irma, paying for their education and continuing to provide for them after their mother's death in 1926. He bought a villa with a garden in Nagytétény, on the outskirts of the city, and had it transferred to the two of them in 1930, while he himself occupied only one room and the veranda; from 1927 the Csáfordi family looked after him. He donated most of his library to the University Library and to the order. He tended the garden himself.\footnote{Obituary, pp.~5\,f. (Suták Papers).}

His connection with the order was not broken, but it was shaped by his life outside the house. In a letter of 17~November 1945 the order's leadership respectfully acknowledged Suták's disposition of his property, while gently reminding him of the needs of the order's young members and praising as an example a donation made in accordance with the order's law---a barely veiled reminder of the vow of poverty.\footnote{Letter from Gyula Zimányi on behalf of the order's leadership to Suták, Budapest, 17~November 1945 (Suták Papers).} Sándor Sík, once his pupil, wrote to him in 1949 as Provincial and in 1951 \enquote{in the name of the whole province}. The order continued to exist, restricted from 1950 to two gymnasia.

His bond with former pupils also remained alive. Graduating classes invited him to their reunions: in 1939 the class of 1894, which called him a \enquote{connecting link} (\emph{összekötő kapocs}) between their schooldays and the present, in 1940 the class of 1895, and in 1951 and as late as June 1954 the class of 1897.\footnote{Letters of 24~June 1939, 28~May 1940, 1951 and 6~June 1954 (Suták Papers). The form of address changes from \emph{Méltóságos} (1939/40) to \emph{Főtisztelendő} (1951/54).} At such reunions, according to the obituary, he spoke in quotations from the Bible; one pupil remembered a sentence from 1934: \enquote{He who finds his faith by way of errors has a faith as firm as a rock.} In old age he no longer came. He drew his pension and thanked God that he received it while so many had lost theirs. Until the end he worked on manuscripts on $n$-dimensional spaces, which according to Pogány consisted almost entirely of formulae; because he had nothing else to give, he handed part of his papers to the Young Pioneers who were collecting waste paper.\footnote{Obituary, pp.~5\,f. (Suták Papers).}

On 22~May 1954 Suták made his will before a notary in Nagytétény; he left his property to Irma Kosaras née Kovács, because she had \enquote{never yet abandoned} him, had always helped him in his old age and, as he hoped, would continue to do so.\footnote{Will of 22~May 1954, certified by Dr.~Ádám Pap (Suták Papers): \enquote{mert [\dots] őt még soha sem hagyta el, öregségében is mindig segítségére volt s reméli, hogy ezután is segítségére lesz.}} The will thus continues a bond of more than thirty years; it is not, as might at first appear, an indictment of the institutions he had served. In the summer of 1954 his strength failed; he fell several times in the garden and in early July spoke confusedly. He died on 19~July 1954, in his eighty-ninth year. The grave he had had prepared ten years earlier in the cemetery of Nagytétény could not be used, because the authorities had closed the cemetery. On 23~July he was buried in Budafok, in the presence of the superior of the house and other members of the order, the Kovács siblings, the parish priest and the neighbours; his former pupils had not learned of the funeral.\footnote{Obituary, p.~6 (Suták Papers).} The grave in the Budafok cemetery (plot~6, row~1, no.~77) survives, though it is now overgrown with ivy. The inscription calls him \enquote{egyetemi tanár} (university professor) and does not mention his membership of the order \citep[221]{ruppert2016}---just as Ottlik had not recognized the Piarist in him. Through the Kosaras family, the part of his papers on which some of the statements in this account rest has been preserved.

\section{Conclusion}\label{sec:conclusion}

Suták's career permits three findings. The first concerns mathematics. From the Bolyai translation and the lecture on geometric axioms (1897/98), through the analysis of the antinomies (1911) and his own axiomatic construction of projective geometry (1914), to the report on von Neumann's axiomatic set theory (1926), one question runs through his work: how the existence and consistency of mathematical objects are to be justified. His answer shifted from the empiricism of the 1898 lecture to a conception of pure geometry as a logical construction. His answers were not original in the sense of new theorems; but they show a Budapest mathematician who knew the international literature on foundations closely, assimilated it independently and conveyed it to a philosophical audience.

The second finding concerns the institution. The appointment of 1911/12 is not the blunder of a provincial faculty but a document of its standards. The majority preferred the proven Privatdozent and teacher trainer to the internationally recognized younger applicants; the ministry additionally appointed the candidate whom the minority opinion had supported with the arguments of solidarity and teacher training. Anyone who wants to understand how twentieth-century Budapest mathematics came into being must also know these standards.

The third finding concerns the political history of the university. In the years of the numerus clausus, Suták, as chairman of the admissions committee, helped to enforce the quota and formulated its criteria; at the same time he defended individual students and his lecture hall against interference. Both follow from the same understanding of service: to his religious order, his school and his university, to order and to rules. Therein lies the paradox of this life. The loyalty that made him a teacher whom his pupils remembered for sixty years also made him a reliable part of an order that excluded others.

Is Suták, then, rightly forgotten? As a mathematician, largely so: he left no results that shaped his field, and even his late work on current questions found no echo. For the history of Budapest mathematics and its university, however, he is not, and this applies precisely to the uncomfortable sides of his work. An account that concentrates on the émigrés loses sight of such figures \citep{frank2009}. Yet they belong to the conditions under which the famous names were trained, examined and appointed---down to the student who in 1921 was probably admitted within the quota proposed by Suták's committee, and whose dissertation bears, in 1926, the signatures of Fejér and Suták.

\printbibliography

\end{document}